\documentclass[notitlepage, twocolumn, letterpaper, prl, nofootinbib, longbibliography, 10pt, showkeys, floatfix]{revtex4-2}

\usepackage{pgfplots}
\usetikzlibrary{intersections}
\pgfplotsset{compat=1.17}

\usepackage{mathtools}   % loads amsmath
\usepackage[normalem]{ulem}
\usepackage{hyperref}
\usepackage[flushleft]{threeparttable}
\usepackage{booktabs}
\usepackage{xcolor}
\usepackage{colortbl} % Needed for gray rows in table
\usepackage{placeins} % allow float barriers

\graphicspath{{figs/}}
\newcommand{\Pmax}{P_\textrm{max}}

\newcommand{\Pstar}{P^{*}}

\newcommand{\Pmem}{P_\textrm{mem}}

\newcommand{\Imax}{I_{\textrm{max}}}
\newcommand{\Pmin}{P_{\textrm{min}}}
\newcommand{\Imin}{I_{\textrm{min}}}
\newcommand{\Ic}{I_{\textrm{c}}}

\newcommand{\Ko}{K_{\textrm{0}}}
\newcommand{\Itot}{I_{\textrm{tot}}}
\newcommand{\Ntot}{N_{\textrm{tot}}}
\newcommand{\Isession}{I^{\textrm{session}}}
\newcommand{\Iactual}{\overline{I}_{\textrm{A}}}
\newcommand{\ICF}{\overline{I}_{\textrm{CF}}}
\newcommand{\TCollab}{T_\textrm{Collab}}
\newcommand{\PCollabNL}{P_{\textrm{Collab}}^{\textrm{NL}}}
\newcommand{\PCollabL}{P_{\textrm{Collab}}^{\textrm{L}}}
\newcommand\ddfrac[2]{\frac{\displaystyle #1}{\displaystyle #2}}

\begin{document}

%************* Front matter ******************
\title{Catalyzing collaborations: Prescribed interactions at conferences determine team formation}

\author{Emma R.~Zajdela}
\email[]{\mbox{emmazajdela@u.northwestern.edu}}
\affiliation{Department of Engineering Sciences and Applied Mathematics, Northwestern University, Evanston, IL, USA}
\author{Kimberly Huynh}
\affiliation{Research Corporation for Science Advancement, Tucson, AZ, USA}
\author{Andy T.~Wen}
\affiliation{Department of Industrial Engineering and Management Sciences, Northwestern University, Evanston, IL, USA}
\author{Andrew L.~Feig}
\affiliation{Research Corporation for Science Advancement, Tucson, AZ, USA}
\author{Richard J.~Wiener}
\affiliation{Research Corporation for Science Advancement, Tucson, AZ, USA}
\author{Daniel M.~Abrams}
\email[]{\mbox{dmabrams@northwestern.edu}}
\affiliation{Department of Engineering Sciences and Applied Mathematics, Northwestern University, Evanston, IL, USA}
\affiliation{Northwestern Institute for Complex Systems, Northwestern University, Evanston, IL, USA}
\affiliation{Department of Physics and Astronomy, Northwestern University, Evanston, IL, USA}

\begin{abstract}
Collaboration plays a key role in knowledge production. Here, we show that patterns of interaction during conferences can be used to predict who will subsequently form a new collaboration, even when interaction is prescribed rather than freely chosen. We introduce a novel longitudinal dataset tracking patterns of interaction among hundreds of scientists during multi-day conferences encompassing different scientific fields over the span of 5 years. We find that participants who formed new collaborations interacted 63\% more on average than those who chose not to form new teams, and that those assigned to a higher interaction scenario had more than an eightfold increase in their odds of collaborating. We propose a simple mathematical framework for the process of team formation that incorporates this observation as well as the effect of memory beyond interaction time. The model accurately reproduces the collaborations formed across all conferences and outperforms seven other candidate models. This work not only suggests that encounters between individuals at conferences play an important role in shaping the future of science, but that these encounters can be designed to better catalyze collaborations.
\end{abstract}


\maketitle

%%%%%%%%%%%%%%%%%%%%%%%%%%%%%%%%%%%%%%
%%%%%%%%% Start main text %%%%%%%%%%%%
%%%%%%%%%%%%%%%%%%%%%%%%%%%%%%%%%%%%%%

The scientific enterprise has increasingly become a team effort \cite{wuchty2007, Fortunatoeaao0185, milojevic2014a, Ahmadpoor13885}. Forming new and interdisciplinary scientific collaborations\footnote{In this paper, we will use the terms ``collaboration,'' ``team formation'' and ``team assembly'' interchangeably, as these terminologies can all be found in the literature designating the same phenomenon \cite{sonnenwald2007}.} is crucial to spurring the innovation necessary to address many current and future challenges facing society \cite{uzzi2013, nsf2018, nih2021}. However, there are intellectual, technical, and logistical obstacles which impede the formation of new teams \cite{Lungeanu2014}.  In particular, research has shown that geographical proximity is a factor in team assembly \cite{wax2017a}. 

Conferences can help overcome these barriers and are one of the main catalysts for the formation of new scientific collaborations. However, convening conferences is expensive in terms of organizational, travel, environmental, and opportunity costs; the direct monetary cost for academic meetings alone is estimated at tens of billions of US dollars each year \cite{klower2020,sarabipour2021}. Given the extensive costs associated with in-person meetings along with benefits virtual conferences may offer in terms of equity and inclusion \cite{skiles2020beyond}, the COVID-19 pandemic has prompted discussions about how scientific conferences should be held even after it is safe to convene them in person \cite{remmel2021}. Some posit that they should continue to be held virtually rather than in person \cite{sarabipour2021}, others that hybrid features should be included \cite{porpiglia2020,margolis2020}, and yet others that the lack of in-person interaction causes a significant damper on scientific productivity and innovation \cite{subramanya2020}. 

Past research has mostly focused on measuring various aspects of the success of scientific collaborations (understood as co-authorship on publications) once formed \cite{funk2017,siudem2020} and the makeup of successful teams (including, e.g., metrics such as the number of institutions present \cite{Jones2008}, team size \cite{wu2019}, and team ``freshness'' \cite{zeng2021}).  There have been some efforts to study scientific team assembly (see, e.g., \cite{milojevic2014a, Guimer2005, Lungeanu2014,freeman2014}), but little is still known about the impact of conferences on collaboration initiation.
Some limited evidence, however, demonstrates that increased interaction among potential members raises the likelihood of team formation \cite{chai2019,campos2018,boudreau2017a}. 

Here we present evidence that properly engineered interaction leads to collaboration, and we go beyond empirical observation by proposing a mathematical model for the origin of this phenomenon.  Such a model has the potential to allow for optimization of conference design to promote collaboration. The model takes as input the pairwise levels of interaction among conference participants as well as their pre-conference familiarity with one other, and estimates the probability that any pair of participants will subsequently form a collaboration.  We test this model using data collected by the Research Corporation for Science Advancement (RCSA) during a series of ``Scialog'' conferences that they organized.

%%%%%%%%%%%%%%%  Subsection  %%%%%%%%%%%%%%%
\section*{Novel ``Scialog'' Dataset}

We construct a novel longitudinal dataset derived from a diverse set of conferences known as ``Scialogs'' \cite{wiener2019}. Organized by the nonprofit funding agency Research Corporation for Science Advancement (RCSA), these conferences seek to accelerate the work of science through research, intensive dialog, community building, and by catalyzing new scientific collaborations on challenges of global significance. Scialog conferences last three days and have an interactive format, with the participation of around 50 fellows, who are invited early-career scientists, and around 10 facilitators, who are more senior scientists. 

For each conference, we have detailed records of how well each participant knew the other participants before the conference, which sessions they attended, with whom they wrote proposals at the end of the conference, and which proposals were funded. Participants are assigned to larger topical discussion sessions of 8--10 fellows facilitated by 1--2 established scientists as well as small group sessions of 3--4 people. Prior knowledge between participants is measured on a four-point scale in a pre-conference survey (see supporting information (SI) for details). At the end of each conference, the participants may self-assemble into teams of 2--4 members to write proposals. A total of 20--35 proposals are typically submitted at each event, of which 5--8 are ultimately funded. 

Scialog conferences are ongoing, so the dataset is still expanding, in particular to include virtual conferences instituted during the pandemic. In its current state, it is comprised of 12 conferences taking place from 2015-2020, divided into four multi-conference series. Each series of conferences deals with a particular scientific topic and most participants return from one conference to the next. For simplicity, in this paper, we have used only data corresponding to the first year of each series to test the model without effects of participants returning in subsequent years. Thus, we used data from four Scialog conferences \cite{wiener2019} from 2015--2018, corresponding to 254 total participants and 4,897 potential pairs of scientists who could form a collaboration\footnote{Some pairings were prohibited because of prior collaborations.}.  For the purpose of this paper, pairs are treated as the fundamental unit of collaboration, and a pair of participants is considered to have formed a collaboration if they were part of the same proposal submission team (triads and tetrads are treated as sets of collaborating pairs). Table \ref{table:ConferenceStatistics} summarizes descriptive statistics for each of the four conferences. Future reports will focus on how the likelihood of collaboration changes across the multi-year arc of each initiative as the individuals get to know one another better, on longer-term collaborations measured by co-authorship, and on effects of virtual versus in-person conferences. 

\begin{table*}[ht]
\centering
\caption{Descriptive Conference Statistics}
\begin{tabular}{crrrrrrrrrr}
Conf. & Year & Topic & Participants & Fellows & \shortstack{Pairs  \\of Fellows} & \shortstack{Proposals \\Submitted} & \shortstack{Proposals \\Funded} & \shortstack{Mean Prior\\Knowledge\\ of Pairs} & \shortstack{ Pairs who\\Collaborated}\\
\midrule
A & 2015 & Molecules Come to Life & 64 & 52 & 1326 & 20 & 5 & 0.91 & 3.5\%\\
B & 2015 & Time Domain Astrophysics & 59 & 49 & 1176 & 30 & 6 & 2.1 & 4.3\%\\
C & 2017 & Advanced Energy Storage & 71 & 60 & 1170 & 35 & 6 & 0.76 & 6.3\%\\
D & 2018 & Chemical Machinery of the Cell & 60 & 50 & 1225 & 24 & 8 & 0.43 & 4.6\%\\
\bottomrule
\label{table:ConferenceStatistics}
\end{tabular}
\vspace{-8mm}
\end{table*}

%%%%%%%%%%%%%%%%%%%%%%%%%%%%%%%%%%%%%%%%%
%%%%%%%%%%%%%%%  Section  %%%%%%%%%%%%%%%
%%%%%%%%%%%%%%%%%%%%%%%%%%%%%%%%%%%%%%%%%
\section*{Empirical Results}

%%%%%%%%%%%%%%%  Subsection  %%%%%%%%%%%%%%%
\subsection*{Interacting More Leads to Collaboration}

\begin{figure}[ht!]
\centering
\includegraphics[width=\linewidth]{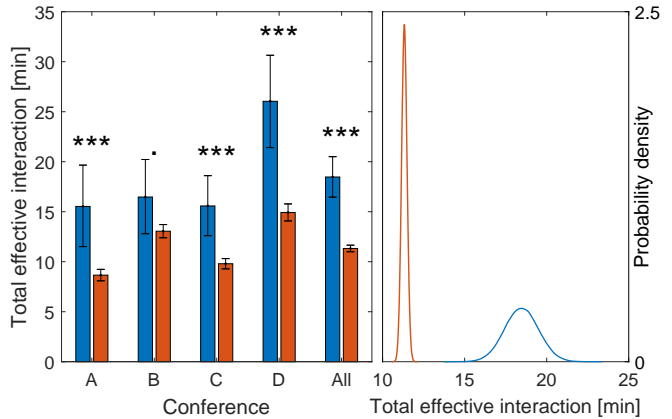}
\caption{\textbf{Effect of interaction on collaboration.} Left panel: blue (left) and red (right) paired bars show bootstrap estimates for mean total effective interaction time for collaborators and non-collaborators, respectively.  Paired bars are shown for each conference analyzed (A-D) as well as the combined data  set of all conferences (All). p-values of the Mann-Whitney U-test for A: $6.0 \times 10^{-4}$, B: $6.3 \times 10^{-2}$, C: $1.4 \times 10^{-4}$, D: $7.7 \times 10^{-6}$, All: $1.0 \times 10^{-12}$. Error bars show mean values of the bootstrapped data with 95\% confidence intervals. Right panel: kernel density estimates showing bootstrapped mean total effective interaction time distributions for collaborators (blue, right) and non-collaborators (red, left) for combined data set of all conferences. 
}
\label{fig:EmpiricalResults}
\end{figure}

We first tested whether pairs who collaborated were different from pairs who did not in terms of total effective interaction $\Itot$ (see appendix~\textbf{\hyperref[sec:DefInteraction]{Defining Interaction}}). To do so, we employed the Mann-Whitney U test \cite{mann1947}; metrics were computed for each conference individually as well as aggregations of the conferences. All metrics indicate that collaborators have significantly more interaction than non-collaborators (see left panel of Fig.~\ref{fig:EmpiricalResults}). The right panel of Fig.~\ref{fig:EmpiricalResults} shows that collaborators spent 63\% more total effective time together than non-collaborators on average across all conferences. This is equivalent to being in a group of 12 people for an extra 45 minutes (60\% of the duration of a topical discussion session) or being in a group of 4 people for an extra 15 min (50\% of the duration of a small group session). 

We also tested whether there was a difference in prior knowledge ($\Ko$) for pairs who collaborated versus those who did not (again using the Mann-Whitney U Test) and found that there was a statistically significant difference at Conference B ($p=6.6 \times 10^{-5}$), Conference C ($p=0.0049$), Conference A at the 10\% level ($p=0.063$) and not at Conference D ($p=0.86$).\footnote{Interestingly, we found no significant effect of interaction or prior knowledge on whether or not a proposal was funded.}

To disentangle causality from correlation for the effect of interaction on collaboration initiation, we performed a test based on counter-factual schedules for one of the conferences. To assign participants to sessions, the conference organizers used a simulated annealing algorithm that attempted to optimize placement based on participant characteristics and information from a pre-conference survey\footnote{The goal was to generate mixtures of participants with varied research methods, disciplines, genders, and minimal professional connections, while also avoiding repeated assignments to groups with similar members in different sessions (thus increasing the number of individuals each participant interacted with).}. Though these assignments were not random, there were numerous options with equivalent or nearly equivalent scores, leaving the organizers to choose among them based on other criteria.  We selected the top 50 alternative solutions for (larger topical) discussion session assignments and the top 50 alternative solutions for small group session assignments, which combined led to 2,500 counter-factual conference schedules.  Note that any one of these schedules could have been chosen as the true conference schedule.

For each counter-factual schedule, we computed the mean interaction\footnote{Here for brevity and clarity we use the word ``interaction'' to mean total effective interaction over the course of the conference.} for all pairs that, in the \textit{actual} conference, ended up collaborating; we denote this as $\ICF^i$ (where integer $i \in [1,2500]$ indexes the particular counter-factual schedule).  It is of interest to compare this quantity to the mean interaction according to the \textit{actual} conference schedule, which we denote  $\Iactual$.  If interacting more had no causal impact on collaboration, we would expect to see little difference in these numbers, with $\ICF^i > \Iactual$ for some $i$ (i.e., under some counter-factual schedules) and $\ICF^i < \Iactual$ for other $i$, in roughly equal proportion.

\begin{figure}[t!]
\centering
\includegraphics[width=\linewidth]{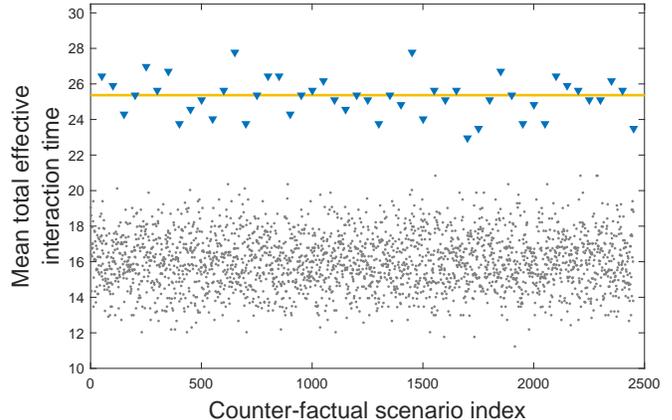}
\caption{\textbf{Mean total effective interaction of collaborators.} Yellow line represents the actual conference, blue inverted triangles are the counter-factual conference solutions where the collaborators were in the same small groups as in the actual conference. Gray points are the other counter-factual conference solutions.}
\label{fig:AlternateUniverse}
\end{figure}

Instead, what we found was that $\Iactual$ was \textit{nearly always} much greater---this was true for more than 99\% of the counter-factual schedules.  The only cases where $\ICF^i > \Iactual$ was observed corresponded to counter-factual scenarios sharing the same exact small-group session assignments but with variations in the larger topical discussion session assignments. Note that this method enabled us to blindly recover the small group assignments knowing only which pairs ultimately collaborated, which strongly suggests a causal connection between intense interaction in a small-group setting and team formation. See Fig.~\ref{fig:AlternateUniverse} for a graphical display of this result.

To quantify the statistical significance of this result, we performed a Wilcoxon signed-rank test \cite{wilcoxon45}. The null hypothesis that the distribution of $\Iactual-\ICF^i$ has zero median is rejected at the $10^{-5}$ level of significance.

\subsection*{Effect Size}

In addition to showing that interaction has a statistically significant effect on collaboration probability, we also wish to know the size of the effect. To evaluate that, we restricted our data to pairs with initial knowledge $\Ko = 0$ (N=984) and used bootstrap statistics to estimate the odds of collaboration for pairs who co-attended one mini session (0.15,  95\% CI [0.10 0.21]) and those who did not co-attend any mini session, but could have in one of the 2,500 counter-factual scenarios (0.017, 95\% CI [0.0085 0.028]). In this case, co-attending a mini-session multiplied the chance of a pair collaborating by 8.7.

%%%%%%%%%%%%%%%%%%%%%%%%%%%%%%%%%%%%%%%%%
%%%%%%%%%%%%%%%  Section  %%%%%%%%%%%%%%%
%%%%%%%%%%%%%%%%%%%%%%%%%%%%%%%%%%%%%%%%%
\section*{Model for the Dynamics of Team Formation at Conferences}

Beyond empirical observations, we develop a mathematical model for the dynamics of team formation at conferences. In 1996, the physicist Serge Galam wrote: ``Do humans behave like atoms?'' \cite{galam1996a}. The mathematical model we present is based on the idea that scientists at a conference behave like molecules in a solution, where formation of a collaboration is analogous to undergoing a chemical reaction. The conference itself acts as a catalyst by lowering the barriers to collaboration and creating more ``productive collisions'' among the scientists who participate. 

\subsection*{Linear Model}As a first simple model, consider a pair of attendees at a conference. We assume that collaboration probability $P(t)$ rises for nonzero interaction intensity $I$ between participants, and when interaction ceases, probability of a collaboration forming decays. For simplicity we assume linear growth and decay processes, leading to the following  ordinary differential equation (ODE) governing the change in collaboration probability over time:
\begin{equation}  \label{eq:model1}
    \frac{dP}{dt} = \underbrace{S \frac{I}{\Imax} \left( 1 - P \right)}_{\textrm{strengthening}} - \underbrace{W P \left(1-\frac{I}{\Imax} \right)}_{\textrm{weakening}}\;.
\end{equation}
This linear model has been constructed to reflect the above assumptions and allow for exponential approach to $P=1$ at ``strengthening'' rate $S$ when $I=\Imax$ (with $\Imax$ the maximum possible interaction intensity), and exponential relaxation to $P=0$ at ``weakening'' rate $W$ when $I=0$. $W$ and $S$ are assumed constant. 

\begin{figure*}[t!]
\centering
\includegraphics[width=11.3cm,height=11.3cm]{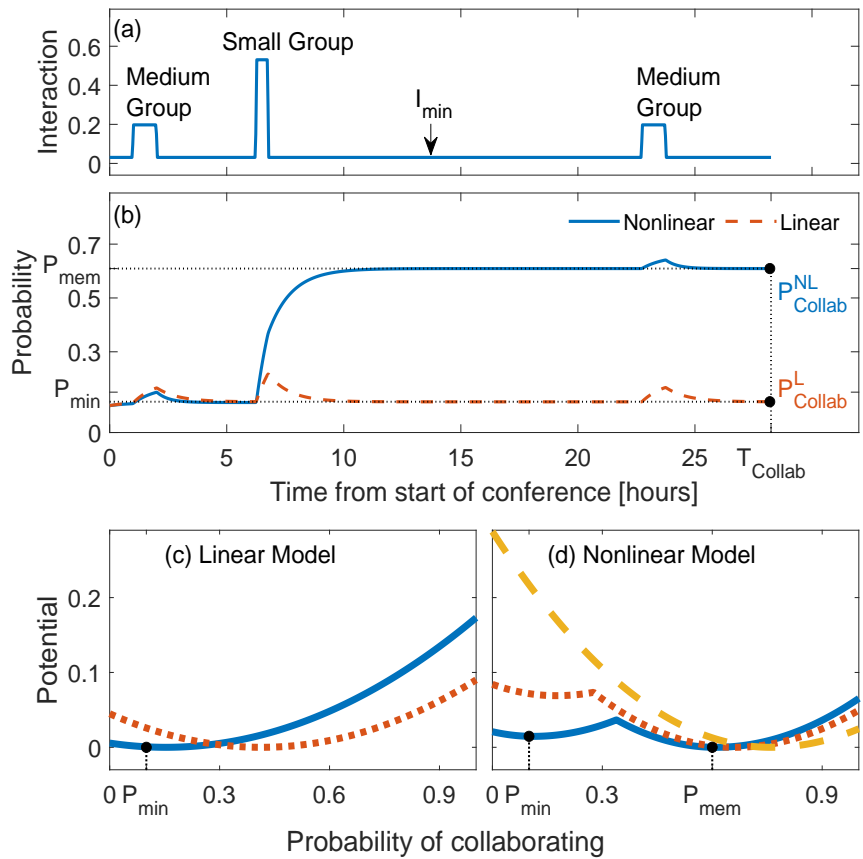}
\caption{\textbf{Model examples.} Effective interaction $I(t)$ (a) and corresponding probability of collaboration (b) as a function of time for a single pair of participants at a realistic conference. Example potential for linear model (c) and nonlinear catalysis model (d). Solid blue curve in (b) shows the probability for the nonlinear catalysis model and dashed red curve shows the probability for the linear model. $\TCollab$ is the time at which teams of participants are formed, $\PCollabNL$ (resp.~$\PCollabL$) is the probability at this time for the nonlinear (resp. linear) model. In panels (c) and (d), solid blue curve shows potential for minimal interaction, dotted red curve shows potential for medium interaction (medium group), dashed yellow curve shows potential for high interaction (small group). Parameter values are the same for linear and nonlinear model where applicable ($a = 0.02, \Ic = 0.2$; $\Imax = 0.6$; $\Pmin = 0.1$; $\Pmem = 0.6$; $\Pmax = 0.9$; $W = 1$; $S = 0.5$.)
}
    \label{fig:realisticconference} 
\end{figure*}

The model can easily be generalized to incorporate more realistic bounds on the collaboration probability:
\begin{equation}  \label{eq:model2}
  \frac {dP}{dt} = S \frac{I}{\Imax} \frac{1 - \frac{P}{\Pmax}}{1 - \frac{\Pmin}{\Pmax}} - W \frac{\left( P-\Pmax \right)}{\Pmax} \left(1 - \frac{I}{\Imax}\right)~.
\end{equation}
Here $\Pmin > 0$ would allow for two participants to form a collaboration even though they have not interacted (for example by being brought together by a third party), and $\Pmax < 1$ allows for non-collaboration even if two participants interact maximally. Two key parameters control the dynamics of this model: the strengthening:weakening ratio and the relative intensity of interaction.

\eqref{eq:model2} can be solved exactly for simple functions $I(t)$.  For example, for constant $I$, $P(t)$ relaxes exponentially to its equilibrium value $\Pstar$. To understand the changes in collaboration probabilities that may occur during the course of a conference, we focus on the case where $I=I(t)$ is not constant, representing the time-varying strength of interaction between two individuals.  Fig.~\ref{fig:realisticconference}(a) shows a ``realistic'' looking example for a 2-day conference, with three sessions of different lengths and intensities, one small (4 people) and two larger (12 people) in this example. Note that $I(t)$ is a dimensionless quantity---for more details on how $I(t)$ was constructed from data, see appendix~\textbf{\hyperref[sec:DefInteraction]{Defining Interaction}} below.

We can express the right hand side of Eq.~\textbf{\ref{eq:model2}} as the derivative of a potential function $V(P)$:
\begin{equation}
    \frac{dP}{dt} = -\frac{\partial V}{\partial P}\;,
    \label{eq:potentialflow}
\end{equation}
where the minus sign implies that stable equilibria occur at potential minima. Since Eq.~\textbf{\ref{eq:model2}} is linear, the resulting potential is quadratic, with a local minimum somewhere in $0 \leq P \leq 1$ (depending on $I$). Fig.~\ref{fig:realisticconference}(c) shows examples of such potentials. 

\vspace{4mm}
\begin{threeparttable}[h]
\centering
\caption{Summary of variables/parameters.}
\begin{tabular}{lll}
\shortstack{Variable or\\Parameter} & Meaning\\
\midrule
$P(t)$ & Probability of collaborating \\
$I(t)$ & Time-varying interaction intensity\\
$\Ko$ & Prior knowledge\\
$a$ & Scaling between interaction and $\Ko$\\
$\Imin$ & Minimum interaction intensity\\
& between sessions\\
\rowcolor{lightgray} 
$\Ic$& Critical interaction intensity at which \\ 
\rowcolor{lightgray} 
& the bifurcation occurs\\
$\Imax$ & Maximum possible interaction intensity\\
$\Itot$ & Total effective interaction\\
$S$  & Rate of strengthening of probability \\
$W$  & Rate of weakening of probability \\
$\Pmin$  & Minimum probability of collaborating \\
\rowcolor{lightgray} 
$\Pmem$  & Memory probability of collaborating \\
$\Pmax$  & Maximum probability of collaborating \\
\bottomrule
\label{table:vars}
\end{tabular}
%  \begin{tablenotes}
%       \small
%       \item This is where authors provide additional information about
%       the data, including whatever notes are needed.
% \end{tablenotes}
\begin{tablenotes}
      \small
      \item Gray rows are parameters only present in the nonlinear catalysis model.
\end{tablenotes}
\end{threeparttable}

%%%%%%%%%%%%%%%  Subsection  %%%%%%%%%%%%%%%
\subsection*{Nonlinear Catalysis Model}

We expect the linear model described above to capture some of the dynamics of formation of new scientific collaborations during conferences. One major limitation, however, is that the linear model relaxes to $\Pmin$ after interaction has ceased, which implies that participants completely forget one other. For a more realistic generalization, we wish to allow scientists who have interacted sufficiently to remember one another long after the interaction has ceased.  To implement this, we modify the potential landscape for a new nonlinear model as shown in Fig.~\ref{fig:realisticconference}(d). 

When interaction $I=0$, there are two stable equilibria, one at the minimum probability $\Pmin$, and the other at memory state $\Pmem$, with an unstable equilibrium in between. As the interaction increases, it acts as a catalyst by changing the shape of the potential function and reducing the barrier between the two stable states. At a critical value of the interaction $I_{c}$ a bifurcation occurs and the barrier disappears, leaving only a single stable equilibrium. If the system gets sufficiently close to that new equilibrium before interaction ceases, the probability will remain permanently in the higher memory state $\Pmem$. 

The exact form of the potential function we employ (a piecewise quadratic function) was motivated by setting rates of strengthening and weakening to be consistent with the original linear model, so curvature of the potential function is higher to the left of each fixed point (strengthening) and lower to the right (weakening). Table \ref{table:vars} summarizes the variables and parameters in the model. See SI for exact form.

%%%%%%%%%%%%%%%%%%%  Subsection  %%%%%%%%%%%%%%%%%%%
\subsection*{Parameter Fitting and Model Selection}

We validate the nonlinear catalysis model by testing how well it explains which pairs of participants ultimately collaborated. The probability of collaborating is the output of the model at time $t=\TCollab$ (see Fig.~\ref{fig:realisticconference}), the start of the period allocated for team formation and proposal writing at the end of the conference. 

\begin{figure}[t!]
\centering
\includegraphics[width=\linewidth]{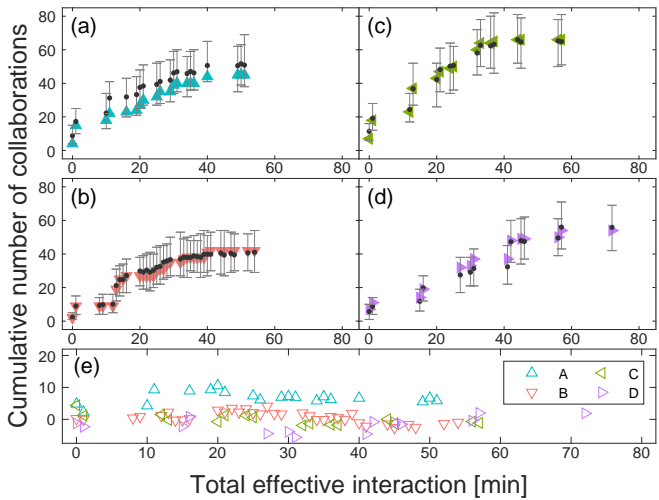}
\caption{\textbf{Test of nonlinear catalysis model.} Sub-figures (a), (b), (c), and (d) show the data (filled triangles) and model predictions (error bars) of the cumulative number of collaborations  as a function of increasing total interaction, measured as total effective interaction at the conference plus a term proportional to prior knowledge. $N=100$ simulations were run for each prediction. Error bars show mean values (points) with 95\% confidence interval. (e) Residuals for model predictions compared to the data.}
\label{fig:ModelData}
\end{figure}

The model was solved numerically with a 4th order Runge-Kutta method assuming initial condition \mbox{$P=\Pmin$} for all pairs. We found optimal parameter values for each conference by minimizing the negative log-likelihood of the model using a constrained Nelder-Mead simplex algorithm \cite{lagarias1998} with a grid of initial values. Fig.~\ref{fig:ModelData} shows agreement between model predictions and data for how collaborations increase with interaction. Of the 85 data points, 84 are contained within the 95\% confidence interval predicted by the model.

We compared the quality of the nonlinear catalysis model to several null models by computing the Akaike Information Criterion (AIC) \cite{akaike1974} and relative likelihood for each one.  These included: (1) a random probability of collaborating between $0$ and $1$ for each pair, (2) a constant probability of collaborating between $0$ and $1$ for all pairs, (3) linear function of prior knowledge $\Ko$, (4) linear function of total effective interaction time $\Itot$, (5) linear function of $\Ko$ and $\Itot$, (6) a threshold-model where the probability was $\Pmem$ if the interaction was ever greater than a critical value $I>\Ic$ and $\Pmin$ if not and (7) the linear model described earlier. Note that candidate models (3)-(7) can be viewed as approximate limits of the nonlinear catalysis model for certain parameter values. 

Table \ref{table:RelativeLikelihood} shows that the AIC of the nonlinear catalysis model is lower compared to the next best model for conferences A, B, C and D, indicating that it is the preferred model for all four conferences.

\begin{table}[t!]
\centering
\caption{Model selection}
\begin{tabular}{crrrr}
Conf. & \shortstack{Nonlinear  \\Model AIC} & \shortstack{Next Best \\ Model AIC} & Next Best Model& \shortstack{Relative\\ Likelihood} \\
\midrule
A & 375.88 & 385.62 & $a \Ko + b \Itot + c$ & 0.0077 \\
B & 337.02 & 340.79 & $a \Ko + b$ & 0.15 \\
C & 526.62 & 541.64 & $a \Ko + b \Itot + c$ & 0.00055 \\
D & 407.96 & 430.53 & $a \Itot + b$ & 1.3e-05\\
\bottomrule
\label{table:RelativeLikelihood}
\end{tabular}
\vspace{-4mm}
\end{table}

%%%%%%%%%%%%%%%%%%%%%%%%%%%%%%%%%%%%%%%%%
%%%%%%%%%%%%%%%  Section  %%%%%%%%%%%%%%%
%%%%%%%%%%%%%%%%%%%%%%%%%%%%%%%%%%%%%%%%%
\section*{Discussion and Limitations}

Our analysis is predicated on the quality of interpersonal interactions as well as their quantity. Impromptu meetings around the coffee maker may differ from those at a conference where participants were encouraged to have a specific conversation and incentivized to form teams though a grant-awarding process, as was the case in the dataset analyzed here. This may limit the contexts in which our model can be applied.

Though we have only tested our model explicitly on the short time scale of one conference, it seems likely that similar dynamics also play out on a multiyear time scale, corresponding to the typical time from project-inception to scientific publication. Future work will address this question explicitly\footnote{Note that some longer time scale memory effects may already be implicitly present in our model through the incorporation of participants' prior knowledge of one another.}.

Our model is probably most applicable when conference participants have limited initial familiarity.  Table \ref{table:ConferenceStatistics} shows that the mean prior knowledge among pairs of participants at conference B was almost five times that of conference D (where participants knew each other the least). That is consistent with the observation that in-person interactions were less significant drivers of collaboration formation at conference B.  In a similar vein, Table \ref{table:RelativeLikelihood} shows that the relative likelihood of the next-best model is highest for conference B, and that next best model is $a \Ko + b$, which accounts only for prior knowledge.

A limitation of this work is that we do not explicitly account for many issues likely affecting team assembly such as, e.g., personality characteristics, homophily, and distance between research fields. However, our approach implicitly incorporates these to some extent through its probabilistic nature. They could also be explicitly incorporated in a more complex future model, but we see the success of the nonlinear catalysis model as remarkable precisely because of its simplicity. 

In our definition of $\Itot$, we have made the simplifying assumption that interaction time is equal for all participants during a group session, despite varying speaking times. Although this approximation does not capture potentially important effects (see, e.g., \cite{mast2002dominance,maclaren2020,khademi2020}), it is a neutral choice absent explicit speaking-time data.

We have also made the simplifying assumption that pairwise dynamics are the primary drivers of collaboration, but more work is needed to quantify the impact of triads, tetrads, or other multilateral interactions.  This may, however, be partially compensated for in our model through the incorporation of parameters $\Pmin$ and $\Imin$.

The Scialog conferences are not representative of all types of conferences, as their goal is explicitly to generate new collaborations and they provide a financial incentive in the form of grants for participants to do so. An extension of this work could study how collaborations are generated at larger conferences or conferences among scientists in regions of conflict where the barrier to collaboration is higher \cite{Zajdela2021}. Research is also needed to compare the benefits and drawbacks of in-person and virtual conferences, especially in terms of the formation of new collaborations. 

\section*{Conclusions}

We have found evidence that prescribed exposure (through structured group interaction) leads to team formation at scientific conferences, and that more interaction leads to a higher probability of forming a collaboration. The nonlinear catalysis model we developed performs better than any other model tested, suggesting that the memory effect it incorporates is key to understanding team formation.  

The nonlinear catalysis model is not necessarily limited to scientific conferences and collaboration---we speculate that it may also have applicable extensions in other areas where matches between individuals within a network are sought. For example, in business settings, employers may wish to promote organic team-building through prescribed sessions among employees. In romantic contexts, a model could inform online dating algorithms and approaches to social interaction. In pedagogical settings, educators might use in-class prescribed group exercises to promote formation of student study groups or teams for collaborative assignments. 

Scientific conferences play an important role in the diffusion of knowledge and can generate novel ideas. We have shown that properly engineered interaction at conferences induces the formation of new scientific collaborations.  Our model helps to illuminate the mechanism by which this occurs, and we hope that it will play a role in designing future conferences such that benefits are maximized and the scientific enterprise is made to proceed as efficiently as possible.

%%%%%%%%%%%%%%%%%%%  Section  %%%%%%%%%%%%%%%%%%%

\appendix*
\section{Appendix: Defining Interaction}
\label{sec:DefInteraction}

For each pair, total effective interaction $\Itot$ is defined as the sum of total effective interaction time during sessions. Total effective interaction during a co-attended session was taken to be proportional to the time one participant spent listening to the other, under the unrealistic but convenient assumption that all participants spoke equally. Thus, for a given pair of participants co-attending a session of time $T_k$ with $N_k$ participants, we assumed: 
\begin{eqnarray}  \label{eq:interactionsession}
    \Isession_k \propto \frac{T_k-T_k/N_k}{N_k-1} \propto \frac{T_k}{N_k}
\end{eqnarray}
where the numerator is the total time spent listening to others and the denominator is the number of people to listen to. When the pair are in different sessions, $\Isession_k=0$. Normalizing so that when $N_k=2$, $\Isession_k=T_k$, we find: $\Isession_k  = \ddfrac{2T_k}{N_k}$. Then with $k$ the index of the session and $m$ the number of sessions (here, $6\leq m\leq8$),

\begin{eqnarray}  \label{eq:Itot}
    \Itot &=& \sum_{k=1}^m \Isession_k  \\
          &=& \sum_{k=1}^m \begin{cases}\ddfrac{2T_k}{N_k} & \textrm{during co-attended session}\nonumber \\
    0 & \textrm{else}
    \end{cases}
\end{eqnarray}

The pairwise interaction intensity profile $I(t)$ was constructed in a similar fashion. To get instantaneous interaction intensity, we divided by the session length $T_k$ and chose units such that a 2-person session would have maximum intensity $\Imax$. We assumed a minimum interaction term (corresponding to informal interaction between sessions) $\Imin = 2/\Ntot$ that depends on the total number of participants $\Ntot$. We then added a term proportional to the initial knowledge $\Ko$\footnote{For simplicity in the model exposition, we considered $\Ko=0$. When $\Ko >0$, the model parameter $\Imax$ needs to be rescaled by the new maximum possible interaction, $6a+1$ (from a hypothetical 2-person session for a pair with maximal $\Ko=6$).}. Eq.~\textbf{\ref{eq:interaction}} summarizes the interaction function as it was implemented. 

\begin{eqnarray}  \label{eq:interaction}
    I(t) = \ddfrac{\Imax}{6a+1} &\Bigg(&a\Ko\\
    &+& \begin{cases} 
           \ddfrac{2}{N}& \textrm{in co-attended sessions} \\
           0 & \textrm{in non co-attended sessions} \\
           \ddfrac{2}{ N_\textrm{tot}} & \textrm{outside of session times} \\
       \end{cases}\Bigg) \nonumber
\end{eqnarray}

Fig.~\ref{fig:realisticconference}(a) shows an example interaction function, with $T=0$ corresponding to one hour before the start of the first topical discussion session (see SI for more example interaction functions).

\begin{acknowledgments}
The authors acknowledge the United States Department of Agriculture NACA 58-3022-0-005 and the Research Corporation for Science Advancement for providing data and assistance. E.R.Z thanks the National Science Foundation Graduate Research Fellowship Program DGE-184216 and the Northwestern Buffett Institute Global Impacts Graduate Research Fellowship for financial support and Maher Said for productive discussions.
\end{acknowledgments}

%\showacknow{} % Display the acknowledgments section

% Bibliography

\bibliography{Zajdela_et_al}

\end{document}

\title{\texorpdfstring{SUPPLEMENTARY INFORMATION FOR\\-----\\ Catalyzing collaborations: Prescribed interactions at conferences determine team formation}{SUPPLEMENTARY INFORMATION FOR: Catalyzing collaborations: Prescribed interactions at conferences determine team formation}}

\author{Emma R.~Zajdela}
\email[]{\mbox{emmazajdela@u.northwestern.edu}}
\affiliation{Department of Engineering Sciences and Applied Mathematics, Northwestern University, Evanston, IL, USA}
\author{Kimberly Huynh}
\affiliation{Research Corporation for Science Advancement, Tucson, AZ, USA}
\author{Andy T.~Wen}
\affiliation{Department of Industrial Engineering and Management Sciences, Northwestern University, Evanston, IL, USA}
\author{Andrew L.~Feig}
\affiliation{Research Corporation for Science Advancement, Tucson, AZ, USA}
\author{Richard J.~Wiener}
\affiliation{Research Corporation for Science Advancement, Tucson, AZ, USA}
\author{Daniel M.~Abrams}
\email[]{\mbox{dmabrams@northwestern.edu}}
\affiliation{Department of Engineering Sciences and Applied Mathematics, Northwestern University, Evanston, IL, USA}
\affiliation{Northwestern Institute for Complex Systems, Northwestern University, Evanston, IL, USA}
\affiliation{Department of Physics and Astronomy, Northwestern University, Evanston, IL, USA}

\begin{abstract}
    \textbf{This file contains:}
    \begin{itemize}
        \item Supplementary text
        \item Figs. 1 to 3
    \end{itemize} 
\end{abstract}

\maketitle

\section{Nonlinear Catalysis Model Potential}

We construct a potential function with the following key properties:
\begin{enumerate}
    \item Two local minima with a barrier between them that lowers, then disappears, as interaction intensity $I$ increases.
    \item When $I=0$, the two local minima should be located at $\Pmin$ and $\Pmem$.
    \item Locally quadratic in $P$ so that the corresponding ODE system ($dP/dt = - \nabla V$) is locally linear (at least in the neighborhood of each minimum).
    \item Asymmetric curvature about each local minimum corresponding to asymmetric strengthening and weakening rate constants $S$ and $W$.
\end{enumerate}

These properties allow for a range of possible potential functions.  We make several additional choices with the goal of producing a simple and tractable function:
\begin{enumerate}
    \item Match form of nonlinear catalysis model with form of linear model for low interaction intensity ($I \ll 1$) and low collaboration probability ($P \ll 1$).
    \item Locate the potential barrier at the midpoint between the two potential minima (i.e., at $P=\half(\Pmin + \Pmem)$) when $I = 0$.
    \item As interaction intensity $I$ increases, deflect both minima linearly and equally in the $+P$ direction such that the upper minimum moves from $\Pmem$ to $\Pmax$ as $I$ goes from $0$ to $\Imax$.
    \item Introduce a tunable parameter $\Ic$ such that the barrier exists for $I<\Ic$ and disappears for $I \ge \Ic$.  To accomplish this, as interaction intensity $I$ increases, deflect the barrier location (local maximum where branches of piecewise function meet) linearly in the $-P$ direction so that the barrier meets the lower potential minimum at interaction intensity $I = \Ic$ (and thus only one equilibrium exists for $I \geq \Ic$).
\end{enumerate}

\begin{figure}[ht!]
    \centering
        \includegraphics[width=0.4\textwidth]{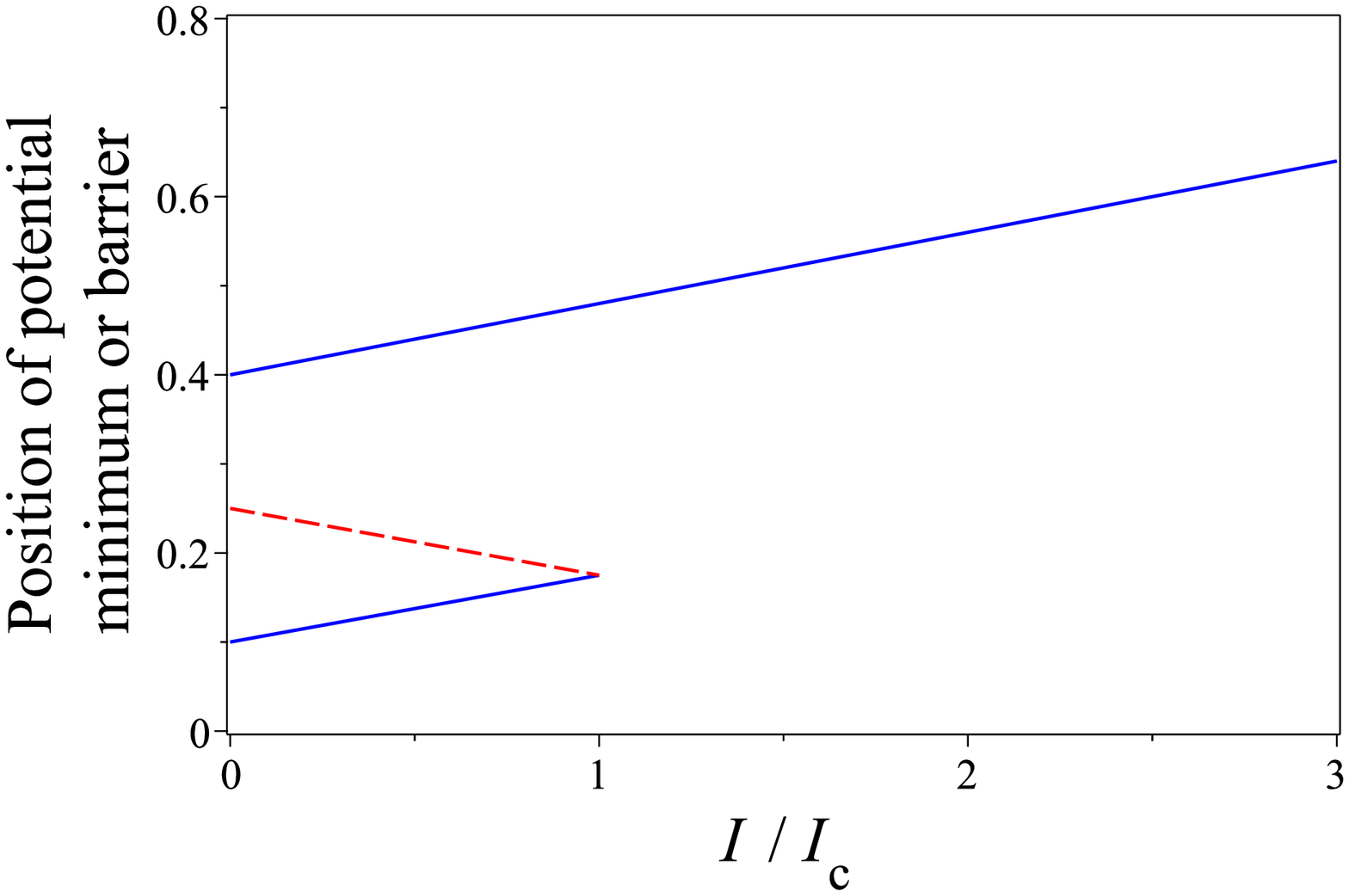}
    \caption{\textbf{Example positions of minima and barrier.} Blue solid lines: Lower and upper potential minima, corresponding to $\Pmin$ and $\Pmem$ when $I=0$.  Red dashed line: Potential barrier position ($\half (\Pmin + \Pmem)$ when $I=0$).  Potential barrier disappears when $I=\Ic$. Constants were chosen for illustration purposes only (here set to $\Pmin=0.1$, $\Pmem=0.4$, $\Pmax=0.8$, $\Ic=1$, $\Imax=5$).}  
    \label{fig:minima}
\end{figure}

\begin{figure}[ht!]
\centering
\includegraphics[width=\textwidth]{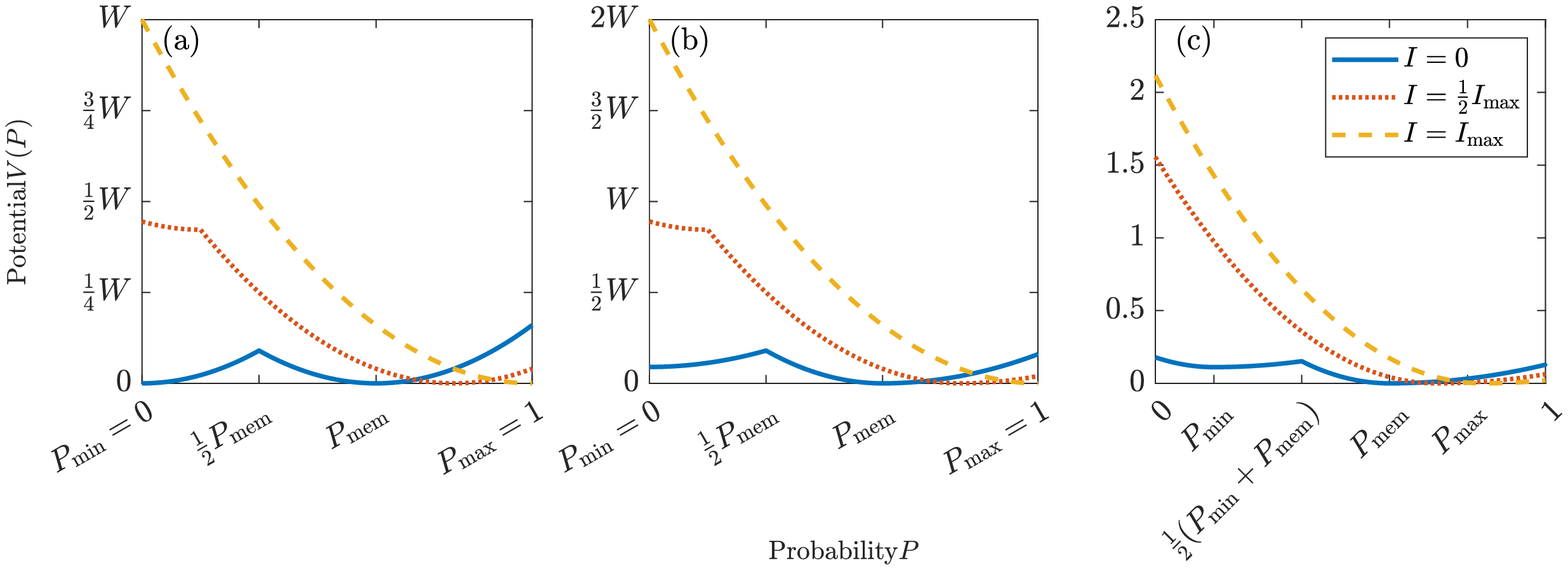}
\caption{\textbf{Potential functions for increasing interaction}. Panel (a) is the simplest version with parameters $\Pmin=0,\Pmax=1,\Imax=1,\Ic=1/2,S=W,W=1,\Pmem=0.6$. Panel (b) is the simplified version with $\Pmin=0, \Pmax=1, \Imax=1, \Ic=1/2, S=2W,W=1,\Pmem=0.6$. Panel (c) is the full version with $\Pmin=0.15, \Pmax=0.8, \Imax=1.2, \Ic=1/2, S=3,W=0.8,\Pmem=0.6$. In each case, the solid blue curve shows the potential for interaction $I=0$, the dotted red curve shows the potential for $I=\half \Imax$ and the dashed yellow curve shows the potential for $I=\Imax$. Note that the vertical axis changes in each panel.}
\label{fig:Potentials}
\end{figure}

See Fig.~\ref{fig:minima} for an illustration of how the potential minima and the potential barrier depend on interaction intensity.  See Fig.~\ref{fig:Potentials} for an illustration of the potential landscapes for varying interaction intensity in the nonlinear catalysis model and each of the below simplifications.

%\end{wrapfigure}
\clearpage

% %%% Each figure should be on its own page

We write the potential function for the nonlinear catalysis model as a piecewise function of $P$ dependent on the interaction intensity $I$.  Although conceptually straightforward, its algebraic representation appears complicated because of its piecewise nature.  Because of that, we first present two simplified special cases to illustrate its structure.

\subsection{Simplest version \texorpdfstring{where $\Pmin=0,\Pmax=1,\Imax=1,\Ic=1/2,S=W$ \\(symmetric strengthening and weakening rates)}{}}\hspace*{\fill} 
\begin{eqnarray}  
    V(P)= \begin{cases} 
    W\left[ (P - \half I \Pmem)^2 - \frac{1}{4} I (\Pmem + 2)(3I\Pmem - 2I - 2\Pmem) \right] & P < \half \Pmem(1 - I) \\
    W\left[ P - I - (1 - I)\Pmem \right]^2 & P \geq \half \Pmem(1 - I) \\
   \end{cases}
\end{eqnarray}
\vspace{5mm}

\subsection{Simplified version \texorpdfstring{where $\Pmin=0,\Pmax=1,\Imax=1,\Ic=1/2,S=2W$ \\(asymmetric strengthening and weakening rates)}{}}\hspace*{\fill} \\
\textbf{First case:} $0 \leq I \leq \Ic$

\begin{equation}
    V_{lowI}(P)= \begin{cases}
        2W\biggl[(P - \half I\Pmem)^2  & P \leq \half I\Pmem \\\hspace{7mm}+ (-\frac{1}{4} \Pmem^2 - \Pmem + 1)I^2+ I\Pmem + \frac{1}{8}\Pmem^2\biggr] \\
    W\biggl[(P - \half I\Pmem)^2 &  \half I\Pmem < P \leq \half (1 - I) \Pmem \\\hspace{7mm}+ (-\half\Pmem^2 - 2\Pmem + 2)I^2 + 2I\Pmem + \frac{1}{4}\Pmem^2\biggr] \\
    2W \left[ P - I - (1 - I)\Pmem  \right]^2 & \half (1 - I) \Pmem < P \leq \Pmem + (1 - \Pmem)I \\
    W \left[ P - I - (1 - I)\Pmem \right]^2 & \Pmem + (1 - \Pmem)I < P
   \end{cases}
\end{equation}

\textbf{Second case:} $\Ic<I\leq\Icint$

\begin{equation} 
    V_{medI}(P)= \begin{cases} 
       2W\left[(P - \half I\Pmem)^2 + (-\frac{1}{4} \Pmem^2 - \Pmem + 1)I^2+ I\Pmem + \frac{1}{8}\Pmem^2\right] & P \leq \Pint \\
      2W \left[ P - I - (1 - I)\Pmem  \right]^2&\Pint <P\leq \Pmem+(1-\Pmem)I \\ W \left[ P - I - (1 - I)\Pmem  \right]^2  & \Pmem+(1-\Pmem)I < P
      
   \end{cases}
\end{equation}

with 

\begin{eqnarray}
    \Pint = -\ddfrac{\Pmem\left[(I^2 - 2I + \frac{7}{8})\Pmem - I^2 + I\right]}{3I\Pmem - 2I - 2\Pmem}
\end{eqnarray}

\textbf{Third case:} $I>\Icint$

\begin{eqnarray} 
    V_{highI}(P)= \begin{cases} 
      2W \left[ P - I - (1 - I)\Pmem  \right]^2& P \leq \Pmem + (1 - \Pmem)I \\
     W \left[ P - I - (1 - I)\Pmem  \right]^2 & \Pmem + (1 - \Pmem)I < P\\
   \end{cases}
\end{eqnarray}
\vspace{5mm}

\subsection{Full version \texorpdfstring{($\Pmin \neq 0$, arbitrary constants)}{}}\hspace*{\fill} \\
\textbf{First case:} $0 \leq I \leq \Ic$

\begin{equation}  
    V_{lowI}(P)= \ddfrac{1}{\Imax^2} \begin{cases} 
      V_1(P) & P\leq \Plow \\
      V_2(P) & \Plow < P \leq \Pmed \\
        S[(\Pmem - P)\Imax + (\Pmax - \Pmem)I]^2 & \Pmed < \Phigh\\
        W[(\Pmem - P)\Imax + (\Pmax - \Pmem)I]^2 & \Phigh < P
   \end{cases}
\end{equation}

with

\begin{eqnarray}
    \Plow = \Pmin+\ddfrac{(\Pmem - \Pmin)I}{4\Ic}\\
    \Pmed = \ddfrac{\Pmin(I + 2\Ic) - \Pmem(I - 2\Ic)}{4\Ic} \\
    \Phigh= \Pmem + \ddfrac{(\Pmax - \Pmem)I}{\Imax} \\
\end{eqnarray}

\begin{eqnarray}
    V_1(P) = \frac{1}{8\Ic^2}\biggl(\{[(8P^2 - 16P\Pmin + 2\Pmem^2 - 4\Pmem\Pmin + 10\Pmin^2)S - 2(\Pmem - \Pmin)^2W]\Ic^2 \nonumber \\ - 4(\Pmem - \Pmin)I[(-\Pmin/2 - \Pmem/2 + P)S - W(\Pmem - \Pmin)]\Ic + I^2(\Pmem - \Pmin)^2(S - 2W)\}\Imax^2 \\+ 4IS\Ic(\Pmem - \Pmin)(\Pmax - \Pmem)(I + 2\Ic)\Imax + 8I^2S\Ic^2(\Pmax - \Pmem)^2\biggr) \nonumber
\end{eqnarray}

and

\begin{eqnarray}
    V_2(P) = \frac{1}{16\Ic^2}\biggl(\{[(4S - 4W)\Pmem^2 - 8\Pmin(S - W)\Pmem + 4S\Pmin^2 + 16(P - \Pmin/2)(P - (3\Pmin)/2)W]\Ic^2 \nonumber\\- 8(\Pmem - \Pmin)((-S/2 - W)\Pmem + PW + \Pmin S/2)I\Ic + I^2(\Pmem - \Pmin)^2(S - 3W)\}\Imax^2 \\+ 8IS\Ic(\Pmem - \Pmin)(\Pmax - \Pmem)(I + 2\Ic)\Imax + 16I^2S\Ic^2(\Pmax - \Pmem)^2\biggr) \nonumber
\end{eqnarray}

\textbf{Second case:} $\Ic<I\leq\Icint$

\begin{equation} 
    V_{medI}(P)= \begin{cases} 
      V_3(P) & P\leq \Pint \\
      \ddfrac{S [(\Pmem-P) \Imax+(\Pmax-\Pmem) I]^2}{\Imax^2} & \Pint < P \leq \Pmem + \ddfrac{(\Pmax - \Pmem)I}{\Imax} \\
        \ddfrac{W[(\Pmem - P)\Imax + (\Pmax - \Pmem)I]^2}{\Imax^2} & \Pmem + \ddfrac{(\Pmax - \Pmem)I}{\Imax} \leq P \\
   \end{cases}
\end{equation}

with 

\begin{equation}
\begin{split}
    \Pint = \ddfrac{1}{16S\Ic \{ [(\Pmax - \Pmem)I + \Imax(\Pmem - \Pmin)]\Ic - \frac{I\Imax(\Pmem - \Pmin)}{4}\})} 
     \{[8S(\Pmem + \Pmin)(\Pmax - \Pmem)I \\+ 6((S + \frac{W}{3})\Pmem + 5(S - \frac{W}{5})\frac{\Pmin}{3})(\Pmem - \Pmin)\Imax]\Ic^2 - 4(\Pmem - \Pmin)[S(\Pmax - \Pmem)I \\+ (((S + 2W)\Pmem + \Pmin(S - 2W))\Imax)/2]I\Ic - I^2\Imax(\Pmem - \Pmin)^2(S - 2W)\} 
\end{split}
\end{equation}

and

\begin{eqnarray}
    V_3(P) = \frac{1}{8\Imax^2\Ic^2} \biggl(\{[(8P^2 - 16P\Pmin + 2\Pmem^2 - 4\Pmem\Pmin + 10\Pmin^2)S - 2(\Pmem - \Pmin)^2W]\Ic^2 \nonumber\\- 4(\Pmem - \Pmin)I((-\frac{\Pmin}{2} - \frac{\Pmem}{2} + P)S - W(\Pmem - \Pmin))\Ic + I^2(\Pmem - \Pmin)^2(S - 2W)\}\Imax^2 \\+ 4IS\Ic(\Pmem - \Pmin)(\Pmax - \Pmem)(I + 2\Ic)\Imax + 8I^2S\Ic^2(\Pmax - \Pmem)^2 \biggr) \nonumber
\end{eqnarray}

\textbf{Third case:} $I>\Icint$

\begin{eqnarray} 
    V_{highI}(P)= \begin{cases} 
      \ddfrac{S[(\Pmem - P)\Imax + (\Pmax - \Pmem)I]^2}{\Imax^2} & P \leq \Pmem + \ddfrac{(\Pmax - \Pmem)I}{\Imax} \\
        \ddfrac{W[(\Pmem - P)\Imax + (\Pmax - \Pmem)I]^2}{\Imax^2 }& \Pmem + \ddfrac{(\Pmax - \Pmem)I}{\Imax} \leq P \\
   \end{cases}
\end{eqnarray}

\textbf{Value of $\Icint$}

$I_1 = I$ such that $\Pint(I) = 0$, $I\in \mathbb{R}$, $0\leq I \leq \Imax$,

$I_2 = I$ such that $\Pint(I) = 1$, $I\in \mathbb{R}$, $0\leq I \leq \Imax$,

$\Icint = \min(I_1,I_2)$.

\newpage
\section{Pre-Conference Survey} \hspace*{\fill} \\

Before attending the conference, participants and fellows were asked to complete the following survey. 

100\% of participants completed the pre-conference survey for conferences A, C and D and 98\% for Conference B.

\begin{quote}

\textbf{Prior Knowledge}

For each name please choose one answer that best describes your relationship with that person prior to this Scialog meeting. There are four categories to choose from: 

\textbf{Unfamiliar:} You are not aware of the research of the person. 

\textbf{Awareness:} Choose this option if you are aware of the research of the person. Examples of "awareness" would be knowing the person's specific area of expertise or knowing details of a recent publication.

\textbf{Discussion:} Choose this option if you have had a substantive discussion about research with this person, through face-to-face conversation, email correspondence, or other means. Please do not select this choice if you have talked with this person and exchanged only basic information about the areas you work in. This level of relationship is meant to be higher than the previous level of "awareness" and presupposes awareness. 

\textbf{Collaborator:} Choose this option if you have ever worked on a project or written a paper together, or formally collaborated with this person on or toward a tangible research output. Please do not select this choice if you have only technically "collaborated" but have never had a substantive research discussion with this person (e.g., coauthored a paper with 100 authors but never interacted). This level of relationship is meant to be higher than the previous level of "discussion" and presupposes awareness and discussion.

Names are listed alphabetically.

Surveys are customized to each respondent. Your name will not appear on the list.

\end{quote}

\begin{quote}

\textbf{Interest in Discussion Topics}

Please choose your interest level for the proposed discussion topics below.
Your input will be used to select the topics for discussion groups at the conference and help us choose which groups you'll be in. Click the "details" button to see more information.

These topics are based on suggestions made by Scialog Fellows, including you, in the conference registration form. Our hope is you will be able to indicate at least a few, and perhaps many, that you are “really into” or would “chime in.”

The order of topics is randomized.
\newline
Respondents are asked to rate each topic on a 5-point scale: No way - Might nap - Would listen - Would chime in - Really into it. 

\end{quote}

\begin{quote}
\textbf{Nominating critical discussion participants}

Listed below are the topics you expressed interest in. For each topic, if you think another Scialog Fellow is an essential person to have in a discussion on that topic, please indicate them below. You may select up to two for each topic but aren't required to select any. Click on the box and start typing or scroll to select a fellow.

\end{quote}

The pre-conference survey results were incorporated into the interaction function as "prior knowledge" $\Ko$ for each pair of fellows (A,B) where $\Ko$ is the sum of prior knowledge reported by A about B and B about A. Thus $\Ko$ for each pair ranges from 0--6 where 0 represents both fellows being unfamiliar with each other and 6 represents both fellows reporting having previously collaborated.

The rules of the Scialog conferences do not allow for participants who have previously collaborated (i.e. pairs with $\Ko \geq 5$) to be on the same proposal submission team. Therefore, when fitting the model to data, we eliminated pairs with $\Ko \geq 5$ (2.1\% of pairs at Conference A, 11.7\% of pairs at Conference B, 3.1\% of pairs at Conference C, 1.5\% of pairs at Conference D).

%\dma{Did we ever restrict to $\Ko = 0$ or $\Ko \leq 3$? Maybe I'm misremembering.}\erz{We restricted to $\Ko <3$ for the version without $\Io$. Log likelihood was lowest for the 4 conferences and AIC was lowest for all but TDA1. We restricted to the case where $\Ko=0$ for what I'm currently working on with effect size but not model fitting. }

\newpage
\section{Group Assignments}\hspace*{\fill} \\

The group assignments were determined prior to the conference with the goal of creating diverse groups for the topical and small group discussion sessions. For the topical discussion groups most of the Fellows who were placed in a group had rated their interest in the topic as a 4 or 5 (on a 1-5 scale with 5 indicating the most interest). Fellows were not placed in a group if they rated the topic under 3. These assignments were accomplished while maintaining diversity in the groups in terms of academic disciplines, research methodologies (e.g. theoretical vs experimental methods), and gender. For the small groups, nearly all Fellows in a group had no previous awareness of the others’ research and none had previously engaged in scientific discussions with the other group members. Participants were mixed so that most small groups included Fellows with different disciplines and methodologies. A simulated annealing algorithm was used to provide candidate groupings based on these criteria. Group assignments for all topical sessions were optimized simultaneously to minimize the same Fellows having repeated assignments together in different sessions. Similarly, all the small group sessions were optimized simultaneously so that no Fellows were ever placed in a small group with another specific Fellow more than once. The algorithm typically returns several solutions with the same or similar energy levels, especially in the case of the small group sessions since they are less constrained. The organizers made the final selection of the group assignments from among the several best solutions. 

Participants who had previously collaborated were not allowed to submit a proposal together, and we therefore eliminated these pairs when fitting models to data. The median percentage of pairs omitted for this reason was 4.0\%.

\newpage
\section{Example interaction patterns}

% %%% Each figure should be on its own page
\begin{figure}[h!]
\centering
\includegraphics[width=\textwidth]{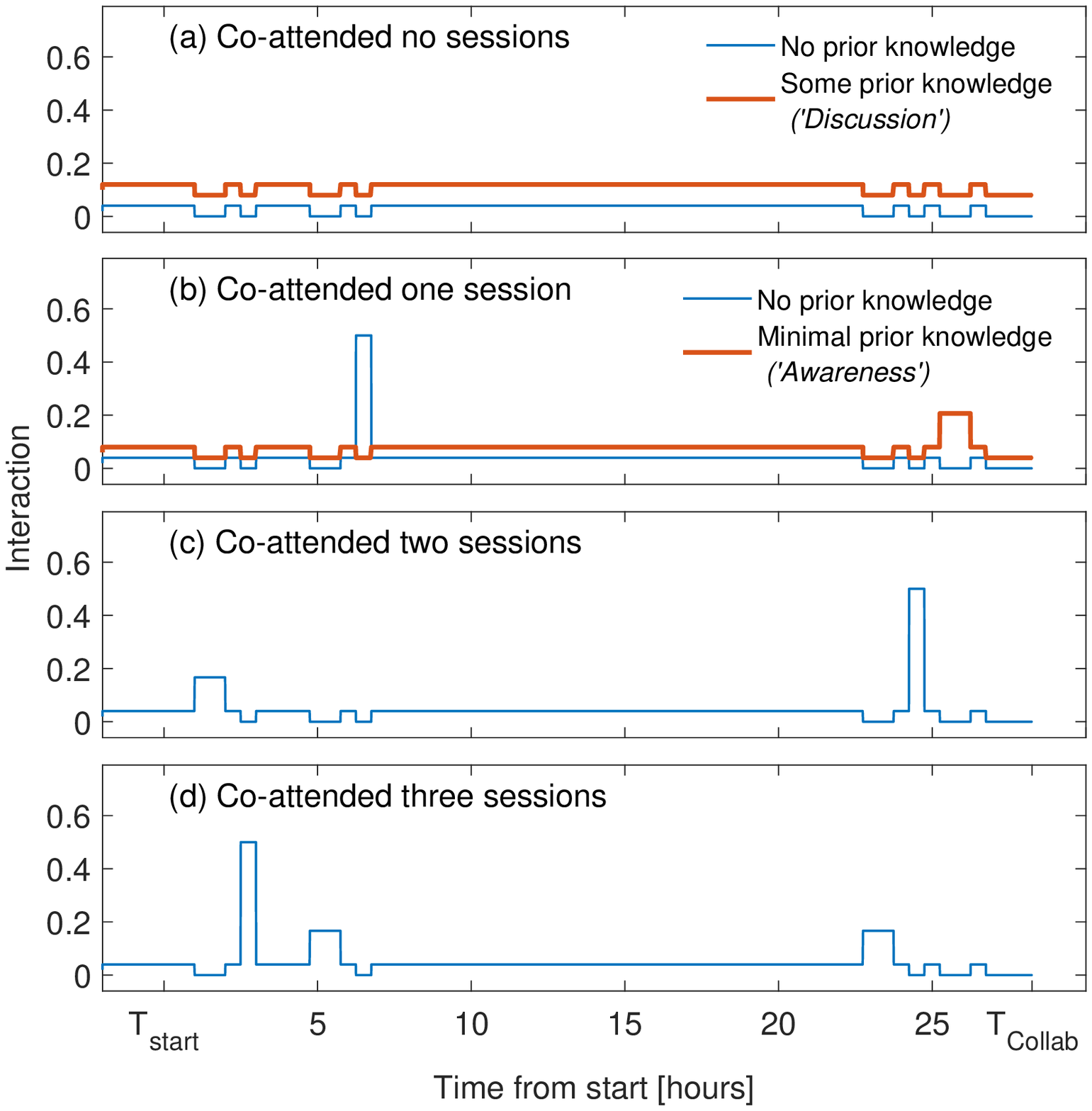}
\caption{\textbf{Example interactions between pairs of participants} from an hour before the first session $T_\textsc{start}$ to the time when proposal writing teams are formed $T_\textsc{collab}$. Higher tophat functions correspond to small group sessions (3-4 people) and medium tophat functions correspond to topical discussion sessions (around 12 people). Drops to zero occur when participants are in different simultaneous sessions.  In panel (b), participants with minimal prior knowledge ($\Ko \leq 2$) may co-attend larger topical sessions but not small group sessions, but in panel (a), participants with some prior knowledge ($\Ko \leq 4$) may not co-attend and sessions. }
\label{fig:Interaction}
\end{figure}